\documentclass[12pt]{article}

\usepackage{amsmath,amssymb,url}
\usepackage{tikz,color,pgf}

\newtheorem{theorem}{Theorem}[section]
\newtheorem{corollary}[theorem]{Corollary}
\newtheorem{lemma}[theorem]{Lemma}
\newtheorem{proposition}[theorem]{Proposition}

\newcommand{\proof}{\noindent{\bf Proof.\ }}
\newcommand{\qed}{\hfill $\square$ \bigskip}
\newcommand{\dstart}{\gamma_g}
\newcommand{\sstart}{\gamma_g^\prime}
\newcommand{\diam}{{\rm diam}}

\title{The domination game played on diameter $2$ graphs}

\author{Csilla Bujt\'as$\/^{a}$, Vesna Ir\v{s}i\v{c}$\/^{a,b}$, Sandi Klav\v{z}ar$^{a,b,c}$, Kexiang Xu$\/^{d}$ \\\\
$^{a}$ \small Faculty of Mathematics and Physics, University of Ljubljana, Slovenia\\
$^{b}$ \small Institute of Mathematics, Physics and Mechanics, Ljubljana, Slovenia \\
$^{c}$ \small Faculty of Natural Sciences and Mathematics, University of Maribor, Slovenia\\
$^{d}$ \small  College of Science, Nanjing University of
 Aeronautics \& Astronautics,\\
 \medskip
  \small Nanjing, Jiangsu 210016, PR China\\
 \small {\tt csilla.bujtas@fmf.uni-lj.si},\quad  \small {\tt vesna.irsic@fmf.uni-lj.si}\\
  \small{\tt sandi.klavzar@fmf.uni-lj.si},\quad   \small {\tt kexxu1221@126.com}
}
\date{\today}

\renewcommand{\gg}{\gamma_{\rm g}}

\begin{document}

\maketitle

\begin{abstract}
Let $\dstart(G)$ be the game domination number of a graph $G$. It is proved that if $\diam(G) = 2$, then $\dstart(G)  \le \left\lceil \frac{n(G)}{2} \right\rceil- \left\lfloor \frac{n(G)}{11}\right\rfloor$. The bound is attained: if $\diam(G) = 2$ and $n(G) \le 10$, then $\dstart(G) = \left\lceil \frac{n(G)}{2} \right\rceil$ if and only if $G$ is one of seven sporadic graphs with $n(G)\le 6$ or the Petersen graph, and there are exactly ten graphs of diameter $2$ and order $11$ that attain the bound.
\end{abstract}


\medskip\noindent
\textbf{AMS Math.\ Subj.\ Class.\ (2020)}: 05C57, 05C69

\section{Introduction}
\label{sec:introduction}

The domination game has been investigated in depth by now, hence let us very quickly recall its definition~\cite{bresar-2010}. The game is played on a graph $G$ by Dominator and Staller who alternately select their vertices. Each selected vertex is required to  dominate at least one new vertex. The game ends when the vertices selected form a dominating set; Dominator's goal is to finish the game as soon as possible, Staller's goal is the opposite. If Dominator is the first to play, we speak of a {\em D-game}, otherwise it is an {\em S-game}. The number of vertices selected in a D-game under the assumption that both players follow optimal strategies is the {\em game domination number} $\dstart(G)$ of $G$.  The corresponding invariant for the S-game is denoted by $\sstart(G)$.

A central theme in the investigation of the game domination number is its upper bounds in terms of the order $n(G)$ of a graph $G$. It all started with the $3/5$-Graph Conjecture~\cite{kinnersley-2013} asserting that if $G$ is an isolate-free graph, then $\dstart(G)\le \frac{3}{5} n(G)$. A strong support for the conjecture is~\cite[Theorem 2.7]{henning-2016}  which asserts that it is true for all graphs with minimum degree at least $2$. The conjecture is still open in general, the best upper bound that holds for all graphs is $\dstart(G)\le \frac{5}{8} n(G)$~\cite[Theorem 2.25]{bujtas-2020a}. Another appealing conjecture, first stated in~\cite{james-2019}, is Rall's $1/2$-Conjecture which asserts that if $G$ is a traceable graph, then $\dstart(G)\le \lceil \frac{1}{2} n(G)\rceil$, cf.~\cite{bujtas-2020b}. For additional topics of interest related to the domination game and its variants see~\cite{borowiecki-2019, bresar-2019, bujtas-2019, dorbec-2015, henning-2017, jiang-2019, nadjafi-2016, ruksasakchai-2019, xu-2018, xu-2018b}.

In this paper we focus on the domination game played on graphs with diameter $2$ and proceed as follows. In the next section we present some preliminary results and introduce a new proof technique (Lemma~\ref{lem:ui}) to be used in the rest of the paper. In Section~\ref{sec:one-half} we prove that if $G$ is a graph with $\diam(G) = 2$, then $\dstart(G)  \le \left\lceil \frac{n(G)}{2} \right\rceil$. Moreover, we show that the equality holds for precisely eight graphs, including the Petersen graph. Based on Section~\ref{sec:one-half}, in the subsequent section we prove that if $G$ is a graph with $\diam(G) = 2$, then $\dstart(G)  \le \left\lceil \frac{n(G)}{2} \right\rceil- \left\lfloor \frac{n(G)}{11}\right\rfloor$. All equality graphs of order $11$ are also discovered. In the concluding section we relate our results to Rall's $1/2$-Conjecture.


\section{Preliminaries}
\label{sec:preliminaries}

We follow the standard graph terminology and notation from~\cite{west-2001}. In particular, if $G$ is a graph, then its  minimum degree, maximum degree, and domination number are denoted by $\delta(G)$, $\Delta(G)$, and $\gamma(G)$, respectively. Also, if $v\in V(G)$, then $N_G(v)$ and $N_G[v]$ denote the open and the closed neighborhood of $v$, respectively. 

First observe that a nontrivial graph $G$ satisfies $\diam(G) \le 2$ if and only if the open neighborhood $N(v)$ is a dominating set for every vertex $v \in V(G)$. Therefore, $\diam(G)=2$ implies $\gamma(G) \leq \delta(G)$. By~\cite[Theorem 1]{bresar-2010}, we have $\dstart(G) \leq 2 \gamma(G) - 1$, and thus $\dstart(G) \leq 2 \delta(G) - 1$. We formulate this observation as a lemma.

\begin{lemma}
	\label{lem:upper-delta}
	If $G$ is a graph of diameter $2$, then $\dstart(G) \leq 2 \delta(G) - 1$.	
\end{lemma}

Recall that if $S \subseteq V(G)$, then $G|S$ denotes a \emph{partially dominated graph}, so a graph $G$, where vertices from $S$ are already dominated. The number of moves remaining in the game on $G|S$ under
optimal play when Dominator, resp.\ Staller, has the next move is denoted by $\dstart(G|S)$, resp.\ $\sstart(G|S)$.

\begin{lemma}
\label{lem:upper-2/3}
If $G = (V, E)$ is a graph of diameter $2$, and $X \subseteq V$ is a non-empty set of (undominated) vertices with $|X|=x$, then the partially dominated graph $G|(V\setminus X)$ satisfies $\dstart(G|(V\setminus X)) \leq \lfloor \frac{2}{3} x + \frac{1}{3} \rfloor$ and $\sstart(G| (V \setminus X)) \leq \lfloor \frac{2}{3} x + \frac{2}{3} \rfloor$.
\end{lemma}

\proof
First, we claim that Dominator can play a vertex which dominates at least two new vertices in each of his moves, except maybe in the last one. Assume that it is Dominator's turn and that he cannot finish the game with a single move. Hence at least two vertices of $G$ are not yet dominated, say $u$ and $v$. Since $\diam(G) = 2$, either $d(u,v) = 1$ or $d(u,v) = 2$. In the first case, Dominator can play $u$ (or $v$) to dominate at least two vertices. And if $d(u,v) = 2$, then Dominator can play a common neighbor of $u$ and $v$, thus again dominating at least two vertices.

Let us now consider the D-game. Assume that the game is played on $G|(V \setminus X)$, and that Dominator uses the above strategy of dominating at least two new vertices at each of his moves (except maybe in his last one), and Staller plays optimally. We distinguish the following two cases.

\medskip\noindent
{\bf Case 1}: The last move of the game is played by Staller. \\
In this case the number of moves played is even, say $2k$, $k\ge 1$. The strategy of Dominator assures that during the game at least $2k + k$ different vertices are dominated. Since in this counting the vertices are pairwise different, we infer that
$$2k + k \le x\,.$$ Since Staller plays optimally, but Dominator maybe not, we can estimate the game domination number as follows:
$$\dstart(G|(V \setminus X)) \le 2k \le \left\lfloor \frac{2}{3} x \right\rfloor \le \left\lfloor \frac{2}{3} x + \frac{1}{3} \right\rfloor\,.$$

\medskip\noindent
{\bf Case 2}: The last move of the game is played by Dominator. \\
Now the number of moves played is odd, say $2k+1$. If $k = 0$, then Dominator's first move finishes the game. Thus $\dstart(G|(V \setminus X)) = 1 \le \lfloor \frac{2}{3} x + \frac{1}{3} \rfloor$. So from now on, we can assume that $k \ge 1$, and hence Dominator cannot finish the game with a single move. By the strategy of Dominator, at least $2k + k + 1$ different vertices are dominated. In this sum, the last $1$ corresponds to the last move of Dominator in which it is possible that he dominates only one new vertex. It follows that $3k + 1 \le x$. Again, as Staller plays optimally but Dominator maybe not, we can estimate that
$$\dstart(G|(V \setminus X)) \le 2k + 1 \le \left\lfloor \frac{2}{3}(x-1) + 1 \right\rfloor =  \left\lfloor \frac{2}{3} x + \frac{1}{3} \right\rfloor\,,$$
and we are done also in this case.

Similar reasoning shows the result for the S-game.
\qed

\begin{corollary}
\label{cor:Delta-bound}
If $G$ is a graph with $\diam(G) = 2$, then
$$\dstart(G)  \le \left\lfloor\frac{2}{3}\Big(n(G) - \Delta(G)\Big)\right\rfloor + 1\,.$$
Moreover, equality holds if $\Delta(G)\in \{n(G)-1, n(G)-2\}$.
\end{corollary}

\proof
 If $\Delta(G) = n(G)-1$, then $\dstart(G) = 1$, and if $\Delta(G)  =  n(G)-2$, then $\dstart(G) = 2$. Hence in both cases the equality holds. In the rest we may thus assume that $\Delta(G)\le n(G)-3$.

Suppose that Dominator starts the game by playing a vertex $v$ of degree $\Delta(G)$. After this move, we are observing an S-game on $G|N[v]$, thus, by Lemma~\ref{lem:upper-2/3}, we have $\sstart(G|N[v]) \leq \lfloor \frac{2}{3}x + \frac{2}{3} \rfloor$, where $x = n(G) - (\Delta(G) + 1)$. This immediately gives $\dstart(G) \leq 1 + \sstart(G|N[v]) \leq 1 + \lfloor \frac{2}{3} (n(G) - \Delta(G)) \rfloor$.
\qed

In~\cite{bresar-2010} it was observed that $\dstart(P) = 5$, where $P$ is the Petersen graph. Hence the equality in Corollary~\ref{cor:Delta-bound} is also sharp for some graphs with $\Delta(G) < n(G) -2$. The following result yields an infinite family of this kind of sharpness examples. Recall that vertices $u$ and $v$ of a graph $G$ are {\em twins} if $N[u] = N[v]$.

\begin{proposition}
\label{prop:Delta-3}
If $G$ is a twin-free graph with $\diam(G) = 2$ and  $\Delta(G) \in \{n(G)-3, n(G)-4\}$, then $\dstart(G) = 3$.
\end{proposition}

\proof
Under the given conditions, Corollary~\ref{cor:Delta-bound} directly implies $\dstart(G) \le 3$. To prove that $\dstart(G) \ge 3$, we need to describe an appropriate strategy of Staller. Assume that Dominator plays a vertex $w$ as his first move. Let $Y = V(G) - N[w]$ and note that $|Y| \ge 2$. If there exist vertices $u,v\in Y$ such that $uv\notin E(G)$, then Staller can play $u$ (or $v$) as her first move, forcing Dominator to play his second move. So assume that $Y$ induces a complete subgraph of $G$ and consider arbitrary vertices $u,v\in Y$. Since $G$ is twin-free, $N(u)\cap N(w) \ne N(v)\cap N(w)$. Let $x\in N(w)$ be a vertex with $xu\in E(G)$ and  $xv\notin E(G)$. If Staller plays $x$ as her first move, then she again  forces Dominator to play one more move. We conclude that $\dstart(G) \ge 3$.
\qed

Let $G_k$, $k\ge 2$, be the graph obtained from the disjoint union of $K_{1,k}$ with leaves $u_1, \ldots, u_k$, and from $K_k$ with vertices $v_1, \ldots, v_k$, by adding the edges $v_iu_j$ for all $i,j\in [k]$, $i\ne j$. (We note in passing that $G_k$ is the Mycielskian $M(K_k)$ of $K_k$~\cite{fisher-1995, mycielski-1955}.) It is straightforward to check that $G_k$  satisfies all the assumptions of Proposition~\ref{prop:Delta-3}, hence constituting an infinite family of equality cases in Proposition~\ref{prop:Delta-3}.

To see that the twin-free condition in Proposition~\ref{prop:Delta-3} is needed, consider the following example. Let $H_k$, $k\ge 4$, be constructed as follows. Start with $K_{k,k}$, and select two adjacent vertices $u$ and $v$ from it. Add a vertex $w$ and connect it to all vertices of $K_{k,k}$. Finally, add two new vertices $x$ and $y$,  the   edge $xy$, and connect both $x$ and $y$ to both $u$ and $v$. Note that $x$ and $y$ are twins, $\diam(H_k) = 2$, and $\Delta(H_k) = n(H_k) - 3$. If Dominator plays $w$ as his first move, then Staller is forced to finish the game with her first move. Hence, $\dstart(H_k) = 2$. To have such examples of $G$ with $\Delta(G) = n(G) - 4$, proceed similarly as above, with the only difference that instead of adding the edge $xy$, we add a triangle.

We conclude this section by proving another lemma, which seems rather technical, but is extremely useful. In a D-game, let $U_i$ and $U_i'$ denote the set of undominated vertices after the move $d_i$ of Dominator and the move $s_i$ of Staller, respectively. Set further $u_i=|U_i|$, $u_i'=|U_i'|$, and $S'_{i}=\{d_1, \dots, d_{i}, s_1, \dots, s_{i}\}.$  A \emph{greedy strategy} of Dominator means that,  for every $i\ge 1$, he plays a vertex $d_i$ that makes $u_i$ as small as possible, that is, Dominator will always select a vertex $v$ that maximizes $|\,N[v] \setminus N[S_i']\,|$.

\begin{lemma}
	\label{lem:ui}
	If $G$ is a graph on $n$ vertices with minimum degree $\delta$ and $i \geq 1$, then
	\begin{align*}
	u_1 & \le n - \delta - 1, \\
	u_{i+1}  & \le  u_i' \left(1-\frac{\delta+1}{n-2i}\right) = u_i' \cdot  \frac{n-2i-\delta - 1}{n-2i},\\
	u_i' & \le u_i - 1,
	\end{align*}
if Dominator follows a greedy strategy.
\end{lemma}

\proof
It is clear that at least $\delta+1$ vertices become dominated with the first move $d_1$ of Dominator, and therefore,
$u_1 \le n-(\delta+1)$. Since Staller must dominate at least one new vertex on each move, we also have $u_i' \le u_i - 1$ for every $i \geq 1$.

It remains to prove the upper bound for $u_{i+1}$. If a vertex $v$ is undominated after the move $s_i$, then $N[v] \subseteq V \setminus S_i'$. As $|N[v]| \geq \delta + 1$, each vertex from $U_i'$ can be dominated by at least $\delta + 1$ different vertices from $V \setminus S_i'$. Since $|V \setminus S_i'| = n - 2i$, we may assume $\delta + 1 \leq n - 2i$ as otherwise the game would be over before the move $d_{i+1}$. A double counting argument shows that $\sum_{v \in U_i'} |N[v]| = \sum_{u \in V\setminus S_i'} |N[u] \cap U_i'|$. As $\sum_{v \in U_i'} |N[v]| \geq u_{i}'(\delta +1)$ and $|V\setminus S_i'|=n-2i$, a vertex from $V\setminus S_i'$ dominates at least $$\frac{u_{i}'(\delta +1)}{n-2i}$$
new vertices on average. Thus,  according to his greedy strategy, Dominator  plays such a vertex that
$$ u_{i+1} \le u'_i - \frac{u_{i}'(\delta +1)}{n-2i} = u_i' \left(1- \frac{\delta+1}{n-2i}\right)\,,$$
which concludes the proof.
\qed

\section{One-half of the order upper bound}
\label{sec:one-half}

In this section we first prove the following result and close the section by discussing an alternative approach to its proof.

\begin{theorem}
\label{thm:1/2-bound}
If $G$ is a graph with $\diam(G) = 2$, then
$$\dstart(G)  \le \left\lceil \frac{n(G)}{2} \right\rceil\,.$$
Moreover, the equality holds if and only if $G$ is one of the graphs from Fig.~\ref{fig:eq_cases_1/2} or the Petersen graph.
\end{theorem}

\begin{figure}[!ht]
	\begin{center}
		\begin{tikzpicture}[thick]
		
		\tikzstyle{every node}=[circle, draw, fill=black!10,
		inner sep=0pt, minimum width=4pt]

        \begin{scope}[yshift = 0cm, xshift = -3.5cm]
		\node (u0) at (0,0) {};
		\node (u1) at (1,0) {};
		\node (u2) at (1.5,1) {};
		\node (u3) at (0.5,2) {};
		\node (u4) at (-0.5,1) {};
		\draw (u0) -- (u1) -- (u2) -- (u3) -- (u4) -- (u0);
		\end{scope}

		\begin{scope}
		\node (u0) at (0,0) {};
		\node (u1) at (1,0) {};
		\node (u2) at (1.5,1) {};
		\node (u3) at (0.5,2) {};
		\node (u4) at (-0.5,1) {};
		\node (u5) at (0.5,1) {};
		
		\draw (u0) -- (u1) -- (u2) -- (u3) -- (u4) -- (u0);
		\draw (u4) -- (u5) -- (u2);
		\end{scope}
		
		\begin{scope}[xshift = 3.5cm]
		\node (u0) at (0,0) {};
		\node (u1) at (1,0) {};
		\node (u2) at (1.5,1) {};
		\node (u3) at (0.5,2) {};
		\node (u4) at (-0.5,1) {};
		\node (u5) at (0.5,1) {};
		
		\draw (u0) -- (u1) -- (u2) -- (u3) -- (u4) -- (u0);
		\draw (u4) -- (u5) -- (u2);
		\draw (u1) -- (u5);
		\end{scope}
		
		\begin{scope}[xshift = 7cm]
		\node (u0) at (0,0) {};
		\node (u1) at (1,0) {};
		\node (u2) at (1.5,1) {};
		\node (u3) at (0.5,2) {};
		\node (u4) at (-0.5,1) {};
		\node (u5) at (0.5,1) {};
		
		\draw (u0) -- (u1) -- (u2) -- (u3) -- (u4) -- (u0);
		\draw (u4) -- (u5) -- (u2);
		\draw (u5) -- (u3);
		\end{scope}

        \begin{scope}[yshift = -2.5cm, xshift = -2cm]
		\node (u0) at (0,0) {};
		\node (u1) at (1.5,0) {};
		\node (u2) at (1.5,1.5) {};
		\node (u3) at (0,1.5) {};
		\draw (u0) -- (u1) -- (u2) -- (u3) -- (u0);
		\end{scope}		
		
        \begin{scope}[yshift = -2.5cm, xshift = 1.3cm]
		\node (u0) at (0,0) {};
		\node (u1) at (1,0) {};
		\node (u2) at (2,0) {};
		\node (v0) at (0,1.5) {};
		\node (v1) at (1,1.5) {};
		\node (v2) at (2,1.5) {};
		\draw (u0) -- (v0) -- (u1) -- (v1) -- (u2) -- (v2) -- (u0);
		\draw (u0) -- (v1);
		\draw (u1) -- (v2);
		\draw (u2) -- (v0);
		
		\end{scope}		

        \begin{scope}[yshift = -2.5cm, xshift = 4.8cm]
		\node (u0) at (0,0) {};
		\node (u1) at (1,0.5) {};
		\node (u2) at (2,0) {};
		\node (v0) at (0,1.5) {};
		\node (v1) at (1,1.0) {};
		\node (v2) at (2,1.5) {};
		\draw (u0) -- (u1) -- (u2) -- (u0);
		\draw (v0) -- (v1) -- (v2) -- (v0);
		\draw (u0) -- (v0);
		\draw (u1) -- (v1);
		\draw (u2) -- (v2);

		\end{scope}

		\end{tikzpicture}
		\caption{Sporadic graphs with $\dstart(G) = \lceil n(G)/2 \rceil$.}
		\label{fig:eq_cases_1/2}
	\end{center}
\end{figure}

\proof
Let $G$ be a graph with $\diam(G) = 2$. Set for this proof $V = V(G)$, $n = n(G)$, $\delta = \delta(G)$, and $\Delta = \Delta(G)$.

Assume first that $\delta \ge n/4 + 1$. By Corollary~\ref{cor:Delta-bound}, keeping in mind that $\delta \le \Delta$, we get:
\begin{align*}
\dstart(G) &  \le \frac{2}{3}\Big(n - \Delta\Big) + 1 \le \frac{2}{3}\Big(n - n/4 - 1\Big) + 1 \\
                & = \frac{n}{2} + \frac{1}{3}\,.
\end{align*}
Since both $n$ and $\dstart(G)$ are integers, we may infer $\dstart(G) \le \lceil \frac{n}{2} \rceil$.

Assume next that $\delta < n/4 + 1$. By Lemma~\ref{lem:upper-delta}, $\dstart(G) < n/2 + 1$. Since both $n$ and $\dstart(G)$ are integers, we get $\dstart(G) \le \lceil \frac{n}{2} \rceil$. This proves the inequality.

For the equality, we have first performed a computer search over all graphs of diameter $2$ and order at most $10$, and found the graphs listed in the statement of the theorem. It thus remains to prove that if $n \ge 11$, then $\dstart(G)  < \left\lceil n/2 \right\rceil$.

First assume that $\delta \le \lfloor \frac{n+1}{4} \rfloor$. Using Lemma~\ref{lem:upper-delta} again we obtain
$$\gamma_g(G) \le 2\delta - 1 \le 2 \left\lfloor \frac{n+1}{4} \right\rfloor -1 < 2\,\frac{n+2}{4}-1 =\frac{n}{2} \le \left\lceil \frac{n}{2} \right\rceil.
$$
For the remaining cases, we now assume $\delta \ge \lfloor \frac{n+1}{4} \rfloor + 1 =  \lfloor \frac{n+5}{4} \rfloor$. Using this bound on the results from Lemma~\ref{lem:ui}, we have that
$$u_1 \le  n-\left\lfloor \frac{n+5}{4} \right\rfloor -1\,,$$
and
$$ u_{i+1} \le (u_i-1)\left(1- \frac{\lfloor \frac{n+9}{4} \rfloor}{n-2i}\right).
$$
Applying these formulas for small cases, we obtain the following conclusions. If $n=11$, then $u_2 \le 2$ and so Staller's move $s_2$ leaves at most one vertex undominated. Therefore, under the greedy strategy of Dominator, the game finishes within $5$ moves that is smaller than $\lceil \frac{11}{2} \rceil$. If $n=12$ or $13$, then $u_2 \le 3$. This gives $u_2' \le 2$. Since $G$ is of diameter $2$, any two vertices are adjacent or share a neighbor, they can be dominated by one move. We conclude $\gamma_g(G) \le 5 < \lceil \frac{12}{2} \rceil < \lceil \frac{13}{2} \rceil$, thus establishing the statement for $n=12$ and $13$. In a similar way we may show $u_2 \le 4$, $u_2' \le 3$ and $u_3 \le 1$  for $n=14$, thus proving $\gamma_g(G) \le 6 < \lceil \frac{14}{2} \rceil$.

From now on, we assume that $n \ge 15$ and, instead of $\delta \ge \lfloor \frac{n+5}{4} \rfloor$, we use the weaker estimation $\delta\ge  \frac{n+2}{4}$. Then, $u_1 \le n-(\delta+1) \le  \frac{3n-6}{4}$ and $u_1' \le \frac{3n-10}{4}$.  By these inequalities and  Lemma~\ref{lem:ui}, we get
\begin{align*}
	u_{2}  & \le\frac{3n-10}{4}  \left(1- \frac{n+6}{4(n-2)}\right) = \frac{(3n-10)(3n-14)}{16(n-2)}.
\end{align*}
Lemma~\ref{lem:upper-2/3} implies that we need at most $1+ 2u_2/3$ further moves to dominate all vertices from $U_2$. Thus, Dominator can ensure that the game finishes in at most
$$4+ \frac{2}{3} \cdot \frac{(3n-10)(3n-14)}{16(n-2)}$$
moves. Now it is enough to consider the strict inequality
$$ 4+ \frac{2}{3} \cdot \frac{(3n-10)(3n-14)}{16(n-2)} < \frac{n}{2}\,.$$
Since it is equivalent to
$$0 < 3n^2-48n +52\,,$$
which is valid for all $n \ge 15$, we may conclude that $\gamma_g (G) < n/2$ holds also for every integer $n \ge 15$.
\qed

One can also think about an alternative approach for proving the upper bound in Theorem~\ref{thm:1/2-bound}  using the following two known results.

\begin{theorem}{\rm \cite[Theorem 3.4]{hellwig-2006}}
\label{thm:old1}
If $G$ is a graph with $\diam(G) = 2$, then $\gamma(G) \leq \lfloor n/4 \rfloor + 1$.
\end{theorem}

\begin{theorem}{\rm \cite[Theorem 5]{meierling-2014}}
\label{thm:old2}
	Let $G$ be a graph of order $n$ and diameter $2$. If $n = 4p+r$ with integers $p \geq 1$ and $0 \leq r \leq 3$, then $\gamma(G) \leq \lfloor n/4 \rfloor = p$, when $r = 0$, $p \geq 4$ or $r = 1$, $p \geq 5$, or $r \in \{2,3\}$, $p \geq 6$.
\end{theorem}

\noindent
Since $\dstart(G) \le 2 \gamma(G) - 1$, see~\cite[Theorem 1]{bresar-2010}, it follows from Theorem~\ref{thm:old1} that for a graph $G$ on $n$ vertices, where $n \not\equiv 0 \pmod 4$, and with diameter $2$, we have $\dstart(G) \leq \lceil n/2 \rceil$. The same conclusion follows from Theorem~\ref{thm:old2} for $n \equiv 0 \pmod 4$ and $n \geq 16$. The remaining cases of graphs on $4$, $8$, and $12$ vertices could then be handled by computer. As stated in the proof of Theorem~\ref{thm:1/2-bound}, we have done this for graphs of order at most $10$. However, the computation for all diameter $2$ graphs on $12$ vertices would require a lot of computer time, hence we did not do it.

\section{A stronger upper bound}
\label{sec:stronger}

In this section we improve Theorem~\ref{thm:1/2-bound} as follows.

\begin{theorem}
	\label{thm:9/22-bound}
	If $G$ is a graph with $\diam(G) = 2$, then
	\begin{equation} \label{eq:2}
	\dstart(G)  \le \left\lceil \frac{n(G)}{2} \right\rceil- \left\lfloor \frac{n(G)}{11}\right\rfloor.
	\end{equation}
	\end{theorem}

\proof
Let $G$ be a graph of diameter $2$ and order $n=n(G)$. First recall that, by Theorem~\ref{thm:1/2-bound}, we have $\dstart(G) \le  \lceil n/2\rceil$ and consequently, (\ref{eq:2}) holds if $n < 11$. By the same theorem, the strict inequality $\dstart(G) <  \lceil n/2\rceil$ is valid whenever $n \ge 11$. The latter implies that  (\ref{eq:2}) is true for every graph $G$ with $11 \le n \le 21$. Therefore, in the following, we assume $n \ge 22$.

The general result~\cite[Theorem 4]{bujtas-2015} directly implies that, in the case of  $\delta(G) \ge 11$, we have $\dstart(G)  < 0.404\, n$. Since
$$ 0.404\, n < \frac{9}{22}\, n = \frac{n}{2} - \frac{n}{11} \le
\left\lceil \frac{n}{2} \right\rceil - \left\lfloor \frac{n}{11}\right\rfloor,$$
we may infer that the theorem holds if $\delta(G) \ge 11$.
\medskip

Assume now that $\delta(G) \le 5$. As $G$ is a graph of diameter~$2$, we have $\dstart(G) \le 2\delta(G)-1 \le 9$ by Lemma~\ref{lem:upper-delta}. On the other hand, $9 \leq 9n/22 \leq \lceil n/2\rceil - \lfloor n/11 \rfloor$ holds under the condition $n \ge 22$, thus proving the theorem for the case of $\delta(G) \le 5$.
\medskip

If $\delta(G) =6$ and $n \ge 27$, we apply the inequalities $\dstart(G) \le 2\delta(G)-1$ and $9n/22 > 11$,  which results in
$$\dstart(G) \le 2\delta(G)-1=11 <  \frac{9}{22}\, n  \le \left\lceil \frac{n}{2} \right\rceil - \left\lfloor \frac{n}{11}\right\rfloor .$$
In addition, one can check that $11 = \lceil n/2 \rceil- \lfloor n/11 \rfloor$ holds, and therefore (\ref{eq:2}) is valid  for $n=25$ and $n=26$.
The remaining cases are $n=22$, $23$ and $24$. Now, Lemma~\ref{lem:ui} moves in. Namely, following the greedy strategy, we observe that Dominator dominates at least $\delta(G)+1=7$ vertices with his first move and therefore, $u_1\le n-7$ and $u_1' \le n-8$. Before the second move of Dominator, there is a vertex that dominates at least
$7u_1'/(n-2)$ vertices. By playing such a vertex, he achieves
$$u_2 \le u_1'\left(1- \frac{7}{n-2}\right) \le \frac{(n-8)(n-9)}{n-2}.$$
After the second move of Staller, we then have $u_2' \le u_2-1$ and may calculate that the number of remaining moves is at most $\lfloor \frac{2}{3}u_2'+ \frac{1}{3} \rfloor$  by Lemma~\ref{lem:upper-2/3}. Thus we conclude
$$\dstart(G) \le 4+ \left\lfloor \frac{2}{3}\left(\frac{(n-8)(n-9)}{n-2}-1\right)+\frac{1}{3}\right\rfloor.$$
The right-hand side formula equals $9$, $10$, $10$, respectively, for $n=22$, $23$, $24$, which are exactly the corresponding values of $\lceil n/2 \rceil- \lfloor n/11 \rfloor$. This finishes the proof for $\delta(G)=6$.
\medskip

Assuming that $\delta(G)=7$, we give a similar reasoning as for $\delta(G)=6$. If $n \ge 32$, then
 $$\dstart(G) \le 2\delta(G)-1=13 <  \frac{9}{22}\, n  \le \left\lceil \frac{n}{2} \right\rceil - \left\lfloor \frac{n}{11}\right\rfloor.$$
 If $29\le n \le 31$, the inequality $13 \le \lceil n/2 \rceil- \lfloor n/11 \rfloor$ is still valid  and so (\ref{eq:2}) is true. Assume finally that $22 \le n \le 28$. It is clear that $u_1 \le n-8$ and $u_1' \le n-9$. Applying Lemma~\ref{lem:ui},  we get
$$u_2 \le u_1'\left(1- \frac{8}{n-2}\right) \le \frac{(n-9)(n-10)}{n-2}.$$
Then, by Lemma~\ref{lem:upper-2/3},
$$\dstart(G) \le 4+ \left\lfloor \frac{2}{3}\left(\frac{(n-9)(n-10)}{n-2}-1\right) + \frac{1}{3}\right\rfloor \le \left\lceil \frac{n}{2} \right\rceil- \left\lfloor \frac{n}{11}\right\rfloor, $$
where the last estimation can be verified by calculating the values for each integer $22 \le n \le 28$. This completes the proof for $\delta(G)=7$.
\medskip

For the next case, suppose $\delta(G)=8$.  If $n \ge 37$, then
$$\dstart(G) \le 2\delta(G)-1=15 <  \frac{9}{22}\, n  \le \left\lceil \frac{n}{2} \right\rceil - \left\lfloor \frac{n}{11}\right\rfloor$$
which proves (\ref{eq:2}). If $n=35$ or $n=36$, then $15 = \lceil n/2 \rceil- \lfloor n/11 \rfloor$ holds  and (\ref{eq:2}) is true. Suppose now that $22 \le n \le 33$. By Lemma~\ref{lem:ui}, $u_1 \le n-9$ and $u_1' \le n-10$ hold. Moreover, Dominator can ensure that
$$u_2 \le u_1'\left(1- \frac{9}{n-2}\right) \le \frac{(n-10)(n-11)}{n-2}\,.$$
Then, taking into account Lemma~\ref{lem:upper-2/3}, we get
$$\dstart(G) \le 4+ \left\lfloor \frac{2}{3}\left(\frac{(n-10)(n-11)}{n-2}-1\right) + \frac{1}{3}\right\rfloor \le \left\lceil \frac{n}{2} \right\rceil- \left\lfloor \frac{n}{11}\right\rfloor$$
by checking the last inequality for each integer between $22$ and $33$. The only remaining case for $\delta(G)=8$ is thus $n=34$. Here, we get $u_2' \le \left\lfloor \frac{24\cdot 23}{32}-1\right\rfloor= 16$ and continue the process with estimating $u_3$  by using Lemma~\ref{lem:ui} again. This yields
$$u_3 \le \left\lfloor u_2'\left(1- \frac{9}{30}\right)\right\rfloor \le \left\lfloor  \frac{16 \cdot 21}{30}\right\rfloor=11$$
and  $u_3' \le 10$. From this point, by Lemma~\ref{lem:upper-2/3},
 Dominator can ensure that the game finishes within seven moves. This establishes $\gamma_g(G) \le 13 <  \lceil 34/2 \rceil- \lfloor 34/11 \rfloor =14$ and completes the proof for $\delta(G)=8$.
\medskip

If $\delta(G)=9$ and $n \ge 42$, we have $\dstart(G) \le 2\delta(G)-1=17 <  \frac{9}{22}\, n $, which implies inequality (\ref{eq:2}). For $n =39$, $40$, $41$, simple calculation shows $17 \le \lceil n/2 \rceil- \lfloor n/11 \rfloor$ and we infer $\dstart(G) \le \lceil n/2 \rceil- \lfloor n/11 \rfloor$ again. To prove the theorem for the remaining cases, $22\le n \le 38$, we first observe that $u_1 \le n-10$ and $u_1' \le n-11$.
Then, by playing greedily, Dominator can ensure (Lemma~\ref{lem:ui})
$$u_2 \le u_1'\left(1- \frac{10}{n-2}\right) \le \frac{(n-11)(n-12)}{n-2}$$
and, by Lemma~\ref{lem:upper-2/3}, we get
$$\dstart(G) \le 4+ \left\lfloor \frac{2}{3}\left(\frac{(n-11)(n-12)}{n-2}-1\right) + \frac{1}{3}\right\rfloor\,.$$
For each integer $22\le n \le 38$, the value of the right-hand side formula is bounded from above by $ \left\lceil \frac{n}{2} \right\rceil- \left\lfloor \frac{n}{11}\right\rfloor$. This completes the proof of (\ref{eq:2}) for $\delta(G)=9$.
\medskip

The last case we have to consider is $\delta(G)=10$. If $n \ge 47$, then $\dstart(G) \le 2\delta(G)-1=19 <  \frac{9}{22}\, n $ holds, thus proving (\ref{eq:2}). For $n=45$ and $n=46$, the inequality $\dstart(G) \le 19 = \lceil n/2 \rceil- \lfloor n/11 \rfloor$
holds, thus implying the statement. If $22\le n \le 44$, we consider the greedy startegy of Dominator which, by Lemma~\ref{lem:ui}, results in $u_1 \le n-11$,  $u_1' \le n-12$,
$$u_2 \le u_1'\left(1- \frac{11}{n-2}\right) \le \frac{(n-12)(n-13)}{n-2}\,.$$
Then, for $22 \le n \le 43$, we can estimate (Lemma~\ref{lem:upper-2/3}) the length of the game as
$$\dstart(G) \le 4+ \left\lfloor \frac{2}{3}\left(\frac{(n-12)(n-13)}{n-2}-1\right) + \frac{1}{3}\right\rfloor \le \left\lceil \frac{n}{2} \right\rceil- \left\lfloor \frac{n}{11}\right\rfloor,$$
 where the last inequality can be easily checked for each integer $n$ in the interval $[22,43]$. If $n=44$, the previous argumentation gives $u_2 \le 32 \cdot 31/42$. Since $u_2$ is an integer, we have $u_2 \le 23$ and $u_2' \le 22$. It follows from Lemma~\ref{lem:ui} that
 $$u_3 \le \left\lfloor u_2'\left(1- \frac{11}{40}\right)\right\rfloor \le \left\lfloor \frac{22\cdot 29}{40} \right\rfloor=15\,,$$
and so, by Lemma~\ref{lem:upper-2/3}, the game will be finished in at most $10$ additional moves. Thus, we conclude $\dstart(G) \le 15$, which implies $\dstart(G) < \lceil 44/2\rceil- \lfloor 44/11 \rfloor =18$.  This completes the proof of Theorem~\ref{thm:9/22-bound}.
\qed

\begin{figure}[ht!]
	\begin{center}
	\includegraphics[trim = 3.8cm 8.8cm 3.8cm 8.5cm, clip, width=.9\textwidth]{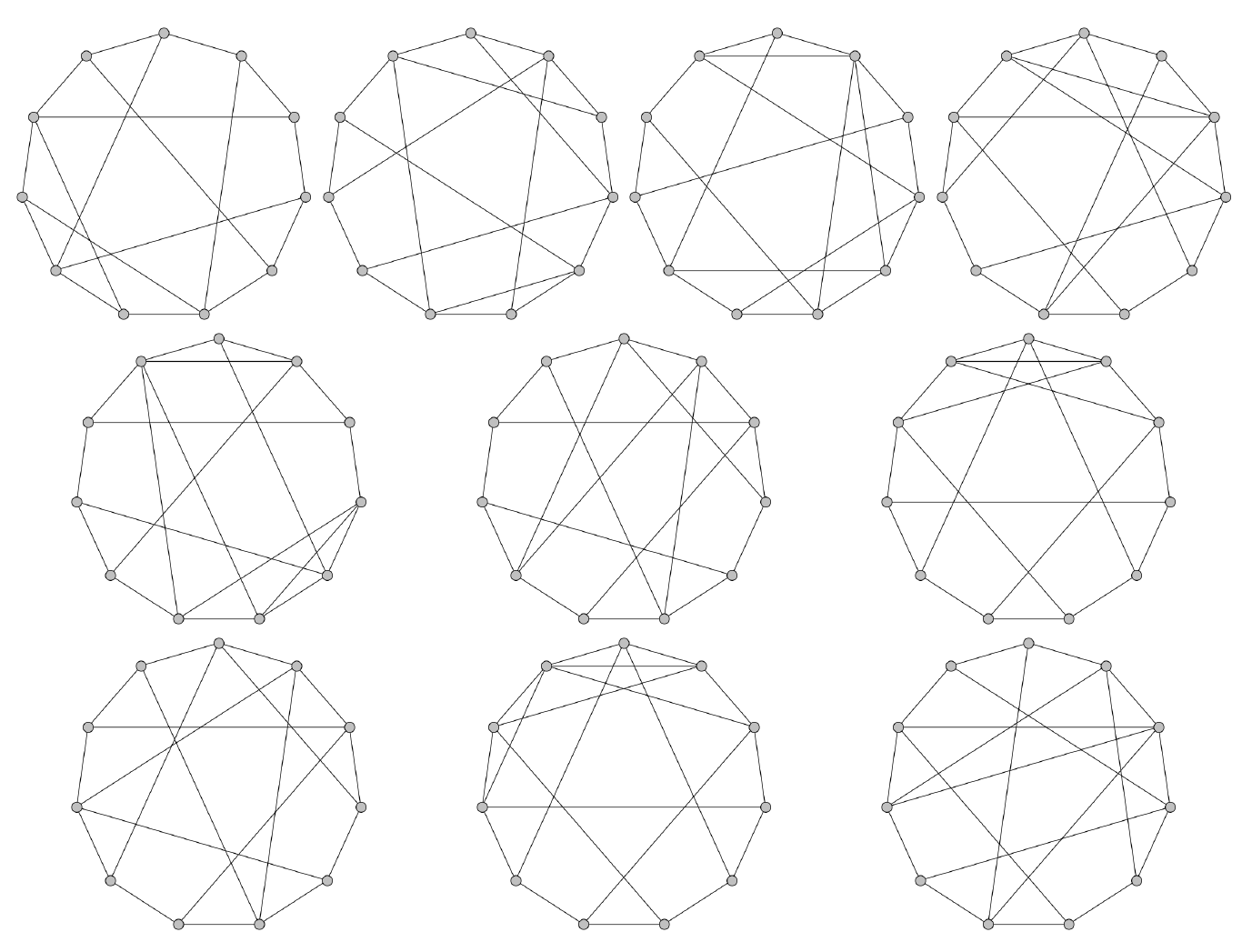}
	\end{center}
	\caption{Graphs on $11$ vertices with diameter $2$ and game domination number $5$.}
	\label{fig:equality11}
\end{figure}

The bound in Theorem~\ref{thm:9/22-bound} is attained. In fact, there are exactly $10$ graphs on $11$ vertices with the game domination number equal to $5$, see Figure~\ref{fig:equality11}. They were obtained using a computer.  Let $G$ be a graph with $n(G) = 11$. If $\delta(G)\ \le 2$, then Lemma~\ref{lem:upper-delta} implies $\dstart(G) \le 3$. And if $\Delta(G) \ge 6$, then Corollary~\ref{cor:Delta-bound} yields $\dstart(G) \le 4$. Hence, in our computer search we only had to check the connected graphs $G$ on $11$ vertices with $\delta(G) \geq 3$ and $\Delta(G) \leq 5$.

\section{Concluding remarks}
\label{sec:conclude}

Note also that the upper bound in Theorem~\ref{thm:1/2-bound} is asymptotically not tight, and for big enough $n$ it follows from know upper bounds on (total) domination number, see for example~\cite{yakoob-2002, desormeaux-2014, dubickas-2021}. In particular, the strongest known result asserts that $\gamma_t(G) < \sqrt{\frac{n \log{n}}{2}} + \sqrt{\frac{n}{2}}$ for all graphs $G$ of diameter $2$ and $n\geq 3$ vertices~\cite[Theorem 1]{dubickas-2021}. Using the well-known bounds $\gg(G) \leq 2 \gamma(G) - 1$ and $\gamma(G) \leq \gamma_t(G)$, this yields 
	$$\gg(G) \leq 2 \left \lfloor \sqrt{\frac{n \log{n}}{2}} + \sqrt{\frac{n}{2}} \right \rfloor - 1.$$
	The latter value is smaller than $\lceil \frac{n}{2} \rceil$ for all $n \geq 65$, and smaller than $\lceil \frac{n}{2} \rceil- \lfloor \frac{n}{11} \rfloor$ for all $n \geq 111$.

Recall that for any fixed positive real number $p < 1$ (which is the probability with which the edges of a random graph are selected mutually independently), almost all graphs are connected with diameter $2$, cf.~\cite[Theorem 13.6]{chartrand-2005}. Hence
Theorem~\ref{thm:9/22-bound} (or Theorem~\ref{thm:1/2-bound} for that matter) imply that
$$\dstart(G) <  \frac{n(G)}{2} $$
holds for almost all graphs $G$.

Theorem~\ref{thm:1/2-bound} and/or Theorem~\ref{thm:9/22-bound} offer another support for Rall's $1/2$-conjecture. That is, the conjecture holds for all graphs with diameter $2$ and consequently for almost all graphs. In this direction, we have tried a different approach than in~\cite{bujtas-2020b}, and proved, using a computer, that the $1/2$-conjecture holds for all Hamiltonian graphs on $n \leq 10$ vertices. Here, we have made use of~\cite[Corollary 5(ii)]{bujtas-2015} to avoid checking graphs with $\delta(G) \geq 5$.

\section*{Acknowledgements}

We are grateful to Ga\v{s}per Ko\v{s}mrlj for providing us with his software that computes game domination invariants. We acknowledge the financial support from the Slovenian Research Agency (research core funding No.\ P1-0297 and projects J1-9109, J1-1693, N1-0095, N1-0108). Kexiang Xu is also supported by NNSF of China (grant No. 11671202) and China-Slovene bilateral grant 12-9.

\end{document}